\documentclass{acm_proc_article-sp}
\newtheorem{theorem}{Theorem}
  \newtheorem{lemma}{Lemma}
 
 \newtheorem{corollary}{Corollary}
 
 \newtheorem{algorithm}{Algorithm}
\begin{document}
\title{An Improved Algorithm for Recovering Exact Value from its
Approximation \titlenote{The work is partially supported by China
973 Project NKBRPC-2004CB318003. }}
\numberofauthors{2}
\author{
\alignauthor Jingzhong  Zhang \\
\affaddr{Laboratory for Automated Reasoning and Programming}\\
\affaddr{Chengdu Institute of Computer Applications}\\
\affaddr{Chinese Academy of Sciences}\\
\affaddr{610041 Chengdu, P. R. China}\\
\email{zjz101@yahoo.com.cn} \alignauthor Yong Feng \\
 \affaddr{Laboratory for Automated Reasoning and Programming}\\
\affaddr{Chengdu Institute of Computer Applications}\\
\affaddr{Chinese Academy of Sciences}\\
\affaddr{610041 Chengdu, P. R. China}\\
\email{yongfeng@casit.ac.cn}}
\maketitle

\begin{abstract}
Numerical approximate computation can solve large and complex
problems fast. It has the advantage of high efficiency. However it
only gives approximate results, whereas  we need exact results in
many fields. There is a gap between approximate computation and
exact results. A bridge overriding the gap was built by Zhang, in
which an exact rational number is recovered from its approximation
by continued fraction method when the error is less than
$1/((2N+2)(N-1)N)$, where $N$ is a bound on absolute value of
denominator of the rational number. In this paper, an improved
algorithm is presented by which a exact rational number is
recovered when the error is less than $1/(4(N-1)N)$.

\end{abstract}
\keywords{Numerical approximate computation, Symbolic-numerical
computation, Continued fraction.}
\section{Introduction}
Numerical approximate computations have the advantage of being
fast, flexible in accuracy and being applicable to large scale
problems. They only give approximate results and are applied in
many fields. However, some fields such as theorem proving, need
exact results and symbolic computations are used to obtain the
exact results. Symbolic computations are principally exact and
stable. They have a high complexity. They are slow and in
practice, are applicable only to small systems.  In recent two
decades, numerical methods are applied in the field of symbolic
computations. In 1985, Kaltofen presented an algorithm for
performing the absolute irreducible factorization, and suggested
to perform his algorithm by floating-point numbers, then the
factor obtained is an approximate one. After then, numerical
methods have been studied to get approximate factors of a
polynomial\cite{corless_C}\cite{Huang}\cite{Mou}\cite{Sasaki_a}
\cite{Sasaki_b}\cite{Sasaki_c}. In the meantime, numerical methods
are applied to get approximate greatest common divisors of
approximate polynomials
\cite{beckermann}\cite{Corless_A}\cite{karmarkar}\cite{corless_D},
to compute functional decompositions\cite{corless_E}, to test
primality\cite{Galligo} and to find zeroes of a
polynomial\cite{Reid}. In 2000, Corless et al. applied numerical
method in implicitization of parametric curves, surfaces and
hypersurfaces\cite{corless_B}. The resulting implicit equation is
still an approximate one.

There is a gap between approximate computations and exact
results\cite{yang}. People usually use rational number
computations to override the gap\cite{Feng}. In fact, these are
not approximate computations but big number computations, which
are also exact computations.  In 2005, Zhang et al proposed an
algorithm to get exact factors of a multivariate polynomial by
approximate computation\cite{zhang} but they did not discuss how
to override the gap.  In \cite{Cheze}, Cheze and Galligo discussed
how to obtain an exact absolute polynomial factorization from its
approximate one, which only involves recovering an integer from
its approximation. They did not discuss how to obtain an exact
rational number from its approximation. Command {\it convert} in
maple can obtain an approximate rational number from a float if we
set variable {\it Digits} to a positive integer. When variable
{\it Digits} is taken to different positive integers, a different
rational numbers are obtained. Which one is the rational number we
want? We do not know. In some cases, what's more, none of the
rational numbers is the one we want. For example, if we want to
get 1/7, we take 0.1196013289 as its approximation. However, no
matter what taking variable {\it Digits} to, we can not obtain 1/7
by command {\it convert}. In \cite{zhang_2}, Zhang and Feng
systematically discussed how to obtain the exact result from its
approximation. They proved that the exact rational number can be
obtained from its approximation when the error is less than
$1/(2N(N-1))$, but in practice, the algorithm requires that the
approximate error is less than $1/((2K+2)N(N-1))$, where $N$ is an
bound on absolute value of denominator of the exact rational
number and $K\ge N$. In this paper, we propose an improved
algorithm which can recover the exact rational number from its
approximation when the error is less than $1/(4N(N-1))$.

The remainder of the paper is organized as follows. Section 2
gives a review of some properties of continued fraction, which are
used to prove our theorems later. Section 3 proves that the
rational number we want is the only one satisfying the error
control $\epsilon\le 1/(2N(N-1))$, and then shows that a kind of
rational numbers can not be obtained from its approximation by
continued fraction method when the error $1/(4N(N-1))<\epsilon\le
1/(2N(N-1))$; and finally, an improved algorithm is proposed.
Section 4 gives some experimental results. The final section makes
conclusions.

\section{Properties of Continued fraction}
A continued fraction representation of a real number is one of the
forms:
\begin{equation}\label{equ:1}
 a_0+\frac{1}{a_1+\frac{1}{a_2+\frac{1}{a_3+\cdots}}},
\end{equation}
where $a_0$ is an integer and $a_1,a_2,a_3,\cdots$ are positive
integers. One can abbreviate the above continued fraction as
$[a_0;a_1,a_2,\cdots]$. For finite continued fractions, note that
$$[a_0;a_1,a_2,\cdots,a_n,1]=[a_0;a_1,a_2,\cdots,a_n+1].$$
So, for every finite continued fraction, there is another finite
continued fraction that represents the same number. Every finite
continued fraction is rational number and every rational number
can be represented in precisely two different ways as a finite
continued fraction. The other representation is one element
shorter, and the final term must be greater than 1 unless there is
only one element. However, every infinite continued fraction is
irrational, and every irrational number can be represented in
precisely one way as an infinite continued fraction. An infinite
continued fraction representation for an irrational number is
mainly useful because its initial segments provide excellent
rational approximations to the number. These rational numbers are
called the convergents of the continued fraction. Even-numbered
convergents are smaller than the original number, while
odd-numbered ones are bigger. If successive convergents are found,
with numerators $h_1,h_2,\cdots$, and denominators
$k_1,k_2,\cdots$, then the relevant recursive relation is:
$$h_n=a_nh_{n-1}+h_{n-2},\,\, k_n=a_nk_{n-1}+k_{n-2}.$$
The successive convergents are given by the formula
$$\frac{h_n}{k_n}=\frac{a_nh_{n-1}+h_{n-2}}{a_nk_{n-1}+k_{n-2}},$$
where $h_{-1}=1$, $h_{-2}=0$, $k_{-1}=0$ and $k_{-2}=1$. Here are
some useful theorems\cite{Fraction}:
\begin{theorem}\label{theo:1}
For any positive $x\in R$, it holds that
\begin{equation}
[a_0,a_1,\cdots,a_{n-1},x]=\frac{xh_{n-1}+h_{n-2}}{xk_{n-1}+k_{n-1}}
\end{equation}
\end{theorem}
\begin{theorem} \label{theo:2}
The convergents of $[a_0,a_1,a_2,\cdots]$ are given by
$$[a_0,a_1,\cdots,a_n]=\frac{h_n}{k_n}$$
and $$k_nh_{n-1}-k_{n-1}h_n=(-1)^n.$$
\end{theorem}
\begin{theorem}\label{theo:add2}
Each convergent is nearer to the $n$th convergent than any of the
preceding convergents. In symbols, if the $n$th convergent is
taken to be $[a_0;a_1,a_2,\cdots,a_n]=x$, then
$$|[a_0;,a_1,a_2,\cdots,a_r]-x|>|[a_0;,a_1,a_2,\cdots,a_s]-x| $$
for all $r<s<n$.
\end{theorem}
{\bf Proof:} Denote $x=[a_0;a_1,a_2,\cdots,a_r,x_{r+1}]$, where
$x_{r+1}=[a_{r+1},\cdots,a_n]$. From theorem \ref{theo:1}, it
holds that
$$x=\frac{x_{r+1}h_r+h_{r-1}}{x_{r+1}k_r+k_{r-1}}$$
Accordingly, we can deduce as follows:
\begin{eqnarray*}
&&x (x_{r+1}k_r+k_{r-1})=x_{r+1}h_r+h_{r-1}\\
&& \Rightarrow x_{r+1}(xk_r-h_r)=-(xk_{r-1}-h_{r-1})
\end{eqnarray*}

Dividing above equation by $x_{r+1}k_r$ yields
$$x-\frac{h_r}{k_r}=(-\frac{k_{r-1}}{x_{r+1}k_r})(x-\frac{h_{r-1}}{k_{r-1}}) $$
Since $x_{r+1}\ge 1$ and $k_r>k_{r-1}>0$, we have
$$0<(\frac{k_{r-1}}{x_{r+1}k_r})<1$$
Therefore, it is proved that
$$|x-\frac{h_r}{k_r}|<|x-\frac{h_{r-1}}{k_{r-1}}| .$$
The proof is finished
\begin{theorem} Let $x=[a_0,a_1,a_2,\cdots]$ and $h_n=a_nh_{n-1}+h_{n-2},\,\,
k_n=a_nk_{n-1}+k_{n-2}$. Then it holds that
$$\left |\frac{h_n}{k_n}-\frac{h_{n-1}}{k_{n-1}}\right |=\frac{1}{k_nk_{n-1}}$$
and
$$\frac{1}{k_n(k_{n+1}+k_n)}<\left |x-\frac{h_n}{k_n}\right|<\frac{1}{k_nk_{n+1}}$$
\end{theorem}

\section{An Improved Algorithm for Recovering the Exact Number from its Approximation}
In this section, we will solve such a problem: Someone has a
positive rational number $\frac{m}{n}$ in his mind, and you only
know an upper bound $N$ of denominator of the rational number, and
he can not tell you the rational number but an approximation of
the rational number at any accuracy. How do you compute the
rational number from one of these approximations? Let's attack
this problem. Without loss of generality, we always assume that
$m$,$n$ are positive numbers. At first, we have a lemma as
follows:
\begin{lemma} \label{lemma:1}
$m,n,p,q$ are integer, and $pn>0$. If
$|\frac{m}{n}-\frac{q}{p}|<\frac{1}{pn}$, then
$\frac{m}{n}=\frac{q}{p}$.
\end{lemma}
{\bf Proof}. $|\frac{m}{n}-\frac{q}{p}|=\frac{|pm-qn|}{pn}$.
Noticing $|pm-qn|$ is a nonnegative integer, and
$|\frac{m}{n}-\frac{q}{p}|<\frac{1}{pn}$, yields $|pm-qn|<1$.
Hence $|pm-qn|=0$. That is $\frac{m}{n}=\frac{q}{p}$. The proof is
finished.
\begin{corollary}\label{cor:1}
$m,n,p,q$ are integers and $p>0$,$n>0$. Let $N\ge\max\{p,n,2\}$.
If $|\frac{m}{n}-\frac{q}{p}|<\frac{1}{N(N-1)}$, then
$\frac{m}{n}=\frac{q}{p}$.
\end{corollary}
{\bf Proof}: When $p\ne n$, it holds that $pn\le N(N-1)$. Hence
$|\frac{m}{n}-\frac{q}{p}|<\frac{1}{N(N-1)}\le \frac{1}{pn}$.
According to lemma \ref{lemma:1}, it is obtained that
$\frac{m}{n}=\frac{q}{p}$. When $p=n$, we have
\begin{eqnarray*}
&&|\frac{m}{n}-\frac{q}{p}|=\frac{|m-q|}{n}<\frac{1}{N(N-1)}\\
&& \Rightarrow |m-q|<\frac{n}{N(N-1)}\le \frac{1}{N-1}\le 1
\end{eqnarray*}
So, it holds that $m=q$. The proof of the corollary is finished.

\begin{theorem}\label{theo:add1}
Let $x=\frac{m}{n}$ be a reduced proper fraction, and
$N\ge\max\{n,2\}$. Assume that $|x-w|<1/(2N(N-1))$. If we get
positive rational number $p/q$ such that $|p/q-w|<1/(2N(N-1))$ ,
where $q\le N$ , then it holds that $x=q/p$.
\end{theorem}
{\bf Proof}: From the assumption of the theorem, we have
$|x-q/p|<1/(N(N-1))$.  According to corollary \ref{cor:1}, it
holds that $q/p=m/n=x$. The proof of the theorem is finished.

Theorem \ref{theo:add1} shows us as follows. One has a rational
number $\frac{m}{n}$ in his mind, if he tells you an approximation
$w$ such that $|\frac{m}{n}-w|<1/(2N(N-1))$, then there is an
unique rational number whose denominator is less than $N$ in the
neighborhood $ (w-1/(2N(N-1)),w+1/(2N(N-1)))$.

The remaining question is as follows. How do we fetch out the
unique rational number in the neighborhood? We wish get it by
continued fraction. Unfortunately, we can not always fetch out the
unique rational number in the neighborhood $
(w-1/(2N(N-1)),w+1/(2N(N-1)))$ by continued fraction method. One
counterexample is the rational number such as $(n-1)/n$ for $n>1$.
Let us show  this: set $N=n$, and its approximation $r :=
(2n+2n^3-4n^2-1)/(2n^2-2n+1)/n$. One can check its error
$d=1/(2n(n-1)+1)<1/(2N(N-1))$. However, one can not recover
rational number $(n-1)/n$ from its approximation by continued
fraction method. In fact. First we can easily compute
$(n-1)/n=[0,1,n-1]$, continued fraction representation of
$n-1)/n$. And then compute continued fraction representation of
$r$ as follows.
\begin{eqnarray*}
&& r=\frac{2n+2n^3-4n^2-1}{(2n^2-2n+1)n}\Rightarrow
\frac{1}{r}=1+\frac{2n^2-n+1}{2n+2n^3-4n^2-1} \\
&& r_1=\frac{2n^2-n+1}{2n+2n^3-4n^2-1}\Rightarrow
\frac{1}{r_1}=n-2+\frac{1+n^2-n}{2n^2-n+1}\\
&& r_2=\frac{1+n^2-n}{2n^2-n+1}\Rightarrow
\frac{1}{r_2}=2+\frac{n-1}{1+n^2-n} \\
&& r_3=\frac{n-1}{1+n^2-n}\Rightarrow
\frac{1}{r_3}=n+\frac{1}{n-1}
\end{eqnarray*}
So, we have $r=[0,1,n-2,2,n,n-1]$. It is obvious that we can not
obtain $n-1/n$ from $[0,1,n-2,2,n,n-1]$. Therefore, we need
smaller neighborhood so as to recover the exact rational number by
continued fraction method.

And now, we discuss how to obtain rational number from its
approximation by continued fraction method. Let $n_2/n_1$ be a
rational number and $r_0$ its approximation. Their continued
fraction representations are $n_2/n_1=[a_0,a_1,\cdots,a_L]$ and
$r_0=[b_0,b_1,\cdots,b_M]$ respectively. We wish that $a_i=b_i$
for $i=0,1,2,\cdots,L-1$ and for the last term of the continued
fraction representations of $n_2/n_1$, either $a_L=b_L$ or
$a_L-1=b_L$, so that we can get $n_2/n_1$ from
$[b_0,b_1,\cdots,b_{L+1}]$. This is the following theorem:
\begin{theorem}\label{theo:3}
Let $n_2/n_1$ be a rational number and $r_0$ its approximation.
Assume that $n_2$,$n_1$ are coprime positive numbers, where
$n_2<n_1$,and $n_1>1$. The representations of $n_2/n_1$ and $r_0$
are $[a_0,a_1,\cdots,a_L]$ and $[b_0,b_1,\cdots,b_M]$
respectively. If $|r_0-n_2/n_1|<1/(4n_1(n_1-1))$, then one of the
following statements must hold.
\begin{itemize}
\item $a_i=b_i$  ($i=0,1,\cdots,L$); \item $a_i=b_i$
($i=0,1,\cdots,L-1$),  $a_L-1=b_L$, and \, $b_{L+1}=1$.
\end{itemize}
\end{theorem}
According to assumption of $n_2<n_1$, we have that $a_0=0$, and
$b_0=0$. Hence $a_0=b_0$. In order to finish the  proof of theorem
\ref{theo:3}, we need two lemmas. Due to $
n_2/n_1=[a_0,a_1,\cdots,a_L]$ and $r_0=[b_0,b_1,\cdots,b_M]$, we
have the following expansions:
\begin{eqnarray}
&&\frac{n_1}{n_2}=a_1+\frac{n_3}{n_2},\;
\frac{n_2}{n_3}=a_2+\frac{n_4}{n_3},\; \cdots , \nonumber \\
&& \frac{n_{L-1}}{n_L}=a_{L-1}+\frac{1}{n_L},\; n_L=a_L
\label{equ:add1}
\end{eqnarray}
and
\begin{eqnarray}
&&\frac{1}{r_0}=b_1+r_1,\; \frac{1}{r_1}=b_2+r_2,\; \cdots
,\nonumber \\
&&\frac{1}{r_{L-1}}=b_L+r_L,\; \cdots,\; \frac{1}{r_{M-1}}=b_M
\label{equ:add2}
\end{eqnarray}
Denoting $d_i=r_i-n_{i+2}/n_{i+1}$, we have a lemma as follows:
\begin{lemma} \label{lemma:2}
Let $n_2/n_1$ be a rational number and $r_0$ its approximation.
Assume that $n_2$,$n_1$ are coprime positive integers, where
$n_2<n_1$,and $n_1>1$. The representations of $n_2/n_1$ and $r_0$
are $[a_0,a_1,\cdots,a_L]$ and $[b_0,b_1,\cdots,b_M]$
respectively. And assume that $a_i=b_i$ for $i\le k<L$($k$ is a
positive integer). Then when $|d_k|<\frac{1}{n_{k+1}(n_{k+1}-1)}$,
it holds that $a_{k+1}=b_{k+1}$ for $k<L-1$; when
$|d_{L-1}|<\frac{1}{n_L(n_L+1)}$, it holds that $a_L=b_L$ or
$a_L-1=b_L$.
\end{lemma}
{\bf Proof}: At first, we show that under the assumption of the
lemma if we have
\begin{equation}\label{equ:add3}
\left|\frac{n_{k+1}^2d_k}{n_{k+2}(n_{k+2}+n_{k+1}d_k)}\right|<
\frac{1}{n_{k+2}}
\end{equation}
then, it holds that $a_{k+1}=b_{k+1}$ for $k<L-1$, and
$a_{k+1}=b_{k+1}$ or $a_{k+1}-1=b_{k+1}$ for $k=L-1$. We discuss
it in two cases:\\ Case 1($k<L-1$): From
$d_k=r_k-n_{k+2}/n_{k+1}$, it holds that
$r_k=n_{k+2}/n_{k+1}+d_k$. Hence we have that
\begin{eqnarray*}
&&\frac{1}{r_k}-\frac{n_{k+1}}{n_{k+2}}=-\frac{n_{k+1}^2d_k}{n_{k+2}(n_{k+2}+n_{k+1}d_k)}\\
\Rightarrow &&
\frac{1}{r_k}=\frac{n_{k+1}}{n_{k+2}}-\frac{n_{k+1}^2d_k}{n_{k+2}(n_{k+2}+n_{k+1}d_k)}\\
\Rightarrow &&
\frac{1}{r_k}=a_{k+1}+\frac{n_{k+3}}{n_{k+2}}-\frac{n_{k+1}^2d_k}{n_{k+2}(n_{k+2}+n_{k+1}d_k)}\\
&&=b_{k+1}+r_{k+1}
\end{eqnarray*}
Hence, obviously,  $a_{k+1}=b_{k+1}$ if and only if
\begin{equation} \label{equ:2} 0\le
\frac{n_{k+3}}{n_{k+2}}-\frac{n_{k+1}^2d_k}{n_{k+2}(n_{k+2}+n_{k+1}d_k)}<1.
\end{equation}
Therefore, if inequality (\ref{equ:add3}) holds, then above
inequality is guaranteed.\\
Case 2:(when $k=L-1$) We have
$$\frac{1}{r_{L-1}}=a_L-\frac{n_{L}^2d_{L-1}}{n_{L+1}(n_{L+1}+n_{L}d_{L-1})}=b_L+r_L $$
From the above equation, if
$$|\frac{n_{L}^2d_{L-1}}{n_{L+1}(n_{L+1}+n_{L}d_{L-1})}|<\frac{1}{n_{L+2}}=1$$
then $a_L=b_L$ for $d_{L-1}<0$, and $a_L-1=b_L$ for $d_{L-1}\ge
0$. Therefore, we have shown that if inequality (\ref{equ:add3})
holds, then $a_{k+1}=b_{k+1}$ for $k<L-1$, and either
$a_{k+1}=b_{k+1}$ or $a_{k+1}-1=b_{k+1}$ for $k=L-1$.

On the other hand, we have
\begin{eqnarray*}
&&\left|\frac{n_{k+1}^2d_k}{n_{k+2}(n_{k+2}+n_{k+1}d_k)}\right|=\frac{n_{k+1}^2|d_k|}{|n_{k+2}(n_{k+2}+n_{k+1}d_k)|}\\
&&=
\frac{n_{k+1}|d_k|}{n_{k+2}(|\frac{n_{k+2}}{n_{k+1}}+d_k|)}\leq
\frac{n_{k+1}|d_k|}{n_{k+2}|(\frac{n_{k+2}}{n_{k+1}}-|d_k|)|}
\end{eqnarray*}
So, in order to ensure inequality (\ref{equ:add3}), we only need
it holds that
\begin{equation}\label{equ:3}
\frac{n_{k+1}|d_k|}{n_{k+2}|(\frac{n_{k+2}}{n_{k+1}}-|d_k|)|}<
\frac{1}{n_{k+2}}
\end{equation}
Solving inequality (\ref{equ:3}) yields
\begin{equation}\label{equ:4}
|d_k|<\frac{n_{k+2}}{n_{k+1}(n_{k+1}+1)}
\end{equation}
When $k<L-1$, we have that $n_{k+2}>1$. So, it holds that
$$\frac{1}{n_{k+1}(n_{k+1}-1)}\leq \frac{n_{k+2}}{n_{k+1}(n_{k+1}+1)} $$
Accordingly, it is obtained that
\begin{equation}
|d_k|<\frac{1}{n_{k+1}(n_{k+1}-1)}
\end{equation}
When $k=L-1$, we have that $n_{L+1}=1$, so it is obtained that
\begin{equation}
|d_{L-1}|<\frac{1}{n_L(n_L+1)}
\end{equation}

The proof of lemma \ref{lemma:2} is finished.

\begin{lemma} \label{lemma:relation_d}
Let $n_2/n_1$ be a rational number and $r_0$ its approximation,
where $n_2$,$n_1$ are coprime positive integers, and $n_2<n_1$,and
$n_1>1$. The continued fraction representations of $n_2/n_1$ and
$r_0$ are $[a_0,a_1,\cdots,a_L]$ and $[b_0,b_1,\cdots,b_M]$
respectively. Denote $d_i=r_i-n_{i+2}/n_{i+1}$ for $i=0,\cdots,L$.
Assume that $a_i=b_i$ for $i\le k<L-1$($k$ is a positive integer
). Then when $|d_k|<\frac{1}{n_{k+1}(n_{k+1}-1)}$, it holds that
\begin{equation}\label{equ:7}
|d_{k+1}|<\frac{n_{k+1}(n_{k+1}-1)}{n_{k+2}(n_{k+2}-1)}|d_k|
\end{equation}
\end{lemma}
{\bf Proof}: Under the assumption that $a_i=b_i$ for
$i=0,1,\cdots,k$, from equation (\ref{equ:add3}), we get
$d_{k+1}=-\frac{n_{k+1}^2d_k}{n_{k+2}(n_{k+2}+n_{k+1}d_k)}$. Hence
we deduce a relation as follows:
\begin{eqnarray*}
|d_{k+1}|&=&\frac{n_{k+1}^2|d_k|}{n_{k+2}^2+n_{k+1}n_{k+2}d_k}=\frac{n_{k+1}|d_k|}{n_{k+2}(\frac{n_{k+2}}{n_{k+1}}+d_k)}\\
&=&\frac{n_{k+1}(n_{k+1}-1)|d_k|}{n_{k+2}(\frac{n_{k+2}(n_{k+1}-1)}{n_{k+1}}+(n_{k+1}-1)d_k)}\\
&=&\frac{n_{k+1}(n_{k+1}-1)|d_k|}{n_{k+2}(n_{k+2}-1+\frac{n_{k+1}-n_{k+2}}{n_{k+1}}+(n_{k+1}-1)d_k)}
\end{eqnarray*}
When $|d_k|<\frac{1}{n_{k+1}(n_{k+1}-1)}$, it holds that
$\frac{n_{k+1}-n_{k+2}}{n_{k+1}}+(n_{k+1}-1)d_k>0$. Hence we have
a relation between $d_{k+1}$ and $d_k$:
\[
|d_{k+1}|<\frac{n_{k+1}(n_{k+1}-1)}{n_{k+2}(n_{k+2}-1)}|d_k| \]
The proof of the lemma is finished.

And now, let us prove the theorem. If
$|d_0|=|r_0-n_2/n_1|<1/(4n_1(n_1-1))$, From lemma
\ref{lemma:relation_d}, we can get
$$|d_i|<1/(4n_{i+1}(n_{i+1}-1))$$ for $i=0,\cdots,L-1$. Note that
$n_L>n_{L+1}=1$ and
$$\frac{1}{4n_{i+1}(n_{i+1}-1)}<\frac{1}{n_{i+1}(n_{i+1}+1)}<\frac{1}{n_{i+1}(n_{i+1}-1)}$$ when
$n_{i+1}>1$. So, it holds that
$$|d_i|<\frac{1}{4n_{i+1}(n_{i+1}-1)}<\frac{1}{n_{i+1}(n_{i+1}+1)}$$
for $i=0,\cdots,L-1$. According to lemma \ref{lemma:2}, the proof
of the theorem is finished.

For an unknown rational number $n_2/n_1$ and its approximation
$r_0$, theorem \ref{theo:3} shows that $n_2/n_1=[b_0,\cdots,b_L]$
or $n_2/n_1=[b_0,b_1,\cdots,b_L,1]$ when
$|r_0-n_2/n_1|<1/(4n_1(n_1-1))$. However, for practical purpose,
we hope the restriction on $n_1>1$ and $n_1>n_2$ can be lifted. So
we have following theorem:
\begin{theorem}\label{theo:practice_alg}
Let $n_0/n_1$ be a reduced rational number and $r$ its
approximation. Assume that $n_0$,$n_1$ are  positive integers and
$N\ge\max\{n_1,2\}$. The continued fraction representations of
$n_0/n_1$ and $r$ are $[a_0,a_1,\cdots,a_L]$ and
$[b_0,b_1,\cdots,b_M]$ respectively. If
$|d|=|r-n_0/n_1|<1/(4N(N-1))$, then one of the following two
statements must hold
\begin{itemize}
\item $a_i=b_i$ for $i=0,\cdots,L$; \item $a_i=b_i$ for
$i=0,\cdots,L-1$, and \, $b_L=a_L-1$, $b_{L+1}=1$.
\end{itemize}
\end{theorem}
{\bf Proof}: We prove the theorem in three cases:\\ Case 1
($n_1>1,\;n_0<n_1$): From  $1/(4N(N-1))\le 1/(4n_1(n_1-1))$ and
theorem \ref{theo:3}, the theorem holds.\\
Case 2($n_1=1$): We have that $a_0=n_0/n_1$. If $d=r-n_0/n_1>0$,
then $b_0=a_0$ and $r_0=r-a_0<1/(4\times 2(2-1))$. If
$d=r-n_0/n_1<0$, then
\begin{eqnarray*}
&&r=a_0-|d|=a_0-1+1-|d|=b_0+1-|d|\\
&&\Rightarrow r_0=r-b_0=1-|d|\\
&&\Rightarrow 1/r_0=1+\frac{|d|}{1-|d|}
\end{eqnarray*}
On the other hand, we have
\begin{eqnarray*}
&&|d|<1/8\Rightarrow 2|d|<1\Rightarrow |d|<1-|d|\\
&&\Rightarrow 0<\frac{|d|}{1-|d|}<1\Rightarrow b_1=1
\end{eqnarray*}

So, we have that $b_0=a_0-1$, $b_1=1$.
\\
Case 3($n_0>n_1>1$): From $n_0/n_1=a_0+n_2/n_1$, it holds that
$n_0/n_1-a_0=n_2/n_1$. On the other hand, we have that
$|n_0/n_1-r|<1/(4N(N-1))\le 1/n_1$. So, we can deduce that
$a_0<r<a_0+1$. Accordingly, it holds that $b_0=a_0$. Hence, we
have
\begin{eqnarray*}
|d|&=&|r-n_0/n_1|=|b_0+r_0-a_0-n_2/n_1|=|r_0-n_2/n_1|\\
&=&d_0<1/(4N(N-1))<1/(4n_1(n_1-1))
\end{eqnarray*}
And now, we have $n_1>1$ and $n_2<n_1$, which is case 1.
Therefore, the proof is finished.

For simplicity, set $s=L$ when
$[a_0,a_1,\cdots,a_L]=[b_0,b_1,\cdots,b_L]$ and $s=L+1$ when
$[a_0,a_1,\cdots,a_L]=[b_0,b_1,\cdots,b_L,1]$. For an unknown
rational number $n_0/n_1$ and its approximation $r$, theorem
\ref{theo:practice_alg} shows that $n_0/n_1=[b_0,\cdots,b_s]$ when
$|r-n_0/n_1|<1/(4N(N-1))$. However, we do not know what the number
$s$ is. The following theorem shows us how to get
$[b_0,\cdots,b_s]$.

\begin{theorem}\label{theo:add3} Let $n_0/n_1$ be a reduced rational number and $r$ its
approximation. Assume that $n_0$,$n_1$ are  positive integers and
$N\ge\max\{n_1,2\}$. We have $r=[b_0,b_1,\cdots,b_M]$ and
$|r-n_0/n_1|<1/(4N(N-1))$. Denote $n_0/n_1=[b_0,b_1,\cdots,b_s]$.
Then, for any positive integer $s<t\le M$, the denominator of
rational number $g=[b_0,b_1,\cdots,b_t]$ is greater than $N$.
\end{theorem}
{\bf Proof:} Proof is given by contradiction. Denote by $m$ the
denominator of $g$. Assume $m\le N$. From theorem \ref{theo:add2},
it holds that $|r-[b_0,b_1,\cdots,b_t]|<|r-[b_0,b_1,\cdots,b_s]|$.
Noting that $|r-[b_0,b_1,\cdots,b_s]|=|r-n_0/n_1|<1/(4N(N-1))$
yields $|r-[b_0,b_1,\cdots,b_t]|<1/(4N(N-1))$. According to
theorem \ref{theo:add1}, it should hold that
$[b_0,b_1,\cdots,b_s]=[b_0,b_1,\cdots,b_t]$. This contradict to
that $t>s$. The proof is finished.

\newcounter{num}
Based on theorem \ref{theo:practice_alg} and theorem
\ref{equ:add3}, an algorithm for obtaining the exact number is as
follows:

\begin{algorithm}\label{alg:1}
Input: a nonnegative floating-point number $r$ and a positive number $N$;\\
Output: a rational number $b$.
\begin{list}{Step \arabic{num}:}{\usecounter{num}\setlength{\rightmargin}{\leftmargin}}
\item Set $i=0$, $tem=r$, $h_{-1}=1$, $h_{-2}=0$, $k_{-1}=0$, and
$k_{-2}=1$; \item \label{step1:2} Get integral part of $tem$ and
assigning it to $a$, assigning its remains to $b$. \item Compute
$h_i=a*h_{i-1}+h_{i-2}$ and $k_i=a*k_{i-1}+k_{i-1}$. If
$k_i>N$,then goto Step \ref{step1:5}; \item Set $i:=i+1$; \item
Set $tem=\frac{1}{b}$ and goto Step \ref{step1:2}; \item
\label{step1:5}  Computing $h_{i-1}/k_{i-1}$ and assigning it to
$b$. \item return $b$.
\end{list}
\end{algorithm}
The correctness of algorithm \ref{alg:1} is obvious from theorem
\ref{theo:practice_alg} and theorem \ref{theo:add3}.

\section{Experimental results}
The following examples  run in the platform of Maple 10 and PIV
3.0G, 512M RAM. They take little time for obtaining exact rational
numbers from their approximations, so we do not show time.

\textbf{Example 1.} Let $a$ be unknown rational number. We only
know a bound of its denominator $N=170$. According to theorem
\ref{theo:practice_alg}, Computing rational number $a$ as follows:
Compute $d=1/(4*N*(N-1))=1/114920$. Assume that we use some
numerical method to get an approximation $b=.8106421859$ such that
$|a-b|<1/d$. Calling algorithm \ref{alg:1} yields as follows.
$$0,
1,\frac{4}{5},\frac{13}{16},\frac{17}{21},\frac{30}{37},\frac{107}{132},\frac{137}{169},\frac{518}{639}$$
When algorithm \ref{alg:1} finds that the denominator of
$\frac{518}{639}$ is greater than $N$, it outputs
$\frac{137}{169}$.

\textbf{Example 2.}Let $a$ be unknown rational number. We only
know a bound of its denominator $N=1790$. According to theorem
\ref{theo:practice_alg}, Computing rational number $a$ as follows:
Compute $d=1/(4*N*(N-1))=1/12809241$. Assume that we use some
numerical method to get an approximation $b=.178870799516605$ such
that $|a-b|<1/d$. Calling algorithm \ref{alg:1} yields as follows:
$$ 0, \frac{1}{5}, \frac{1}{6}, \frac{2}{11}, \frac{5}{28}, \frac{17}{95},
\frac{22}{123}, \frac{149}{833}, \frac{171}{956},
\frac{320}{1789}, \frac{1131}{6323}$$ When the algorithm finds
that the denominator of $1131/6323$ is greater than $N$, it
outputs $320/1789$.

\textbf{Example 3.}Let $a$ be unknown rational number. We only
know a bound of its denominator $N=18$. According to theorem
\ref{theo:practice_alg}, Computing rational number $a$ as follows:
Compute $d=1/(4*N*(N-1))=1/1225$. Assume that we use some
numerical method to get an approximation $b=1.881536615$ such that
$|a-b|<1/d$. Calling algorithm \ref{alg:1} yields as follows.
$$1, 2, 15/8, 32/17, 111/59$$
When the algorithm finds that the denominator of $111/59$ is
greater than $N$, $32/17$.

\textbf{Example 4.}  This example is an application in obtaining
exact factors from their approximations. Let
$p=-16-56y-48z+64x^2-32xy+48xz-45y^2-96yz-27z^2$ be a polynomial.
We want to use approximate method to get its exact factors over
rational number field. First, we transform $p$ to a monic
polynomial as follows:
$$ p=x^2-\frac{1}{2}xy+\frac{3}{4}xz-\frac{45}{64}y^2-\frac{3}{2}yz-\frac{27}{64}z^2-\frac{7}{8}y-\frac{3}{4}z-\frac{1}{4}$$
the least common multiple of denominators of coefficients of
polynomial $p(x,y,z)$ is 64, which is an upper bound\cite{zhang_1}
of denominators of coefficients of the monic factors of polynomial
$p$. Taking $N=65$ yields $d=1/(4*65*64)=1/16128$. We use
numerical methods to get its approximate factors as
follows\cite{zhang_1}:
$$\bar g_1=1.0000x+.6250000000067y+1.124999999530z+.50000$$
$$\bar g_2=1.0000x-1.125000000015y-.3749999995480z-.50000$$
the error of coefficients of $\bar g_1$ and $\bar g_2$ is less
than $d$ by the numerical methods. According to theorem
\ref{theo:practice_alg}, taking $N$ in algorithm \ref{alg:1}, we
obtain two exact factors:
$$g_1=x+\frac{5}{8}y+\frac{9}{8}z+\frac{1}{2}$$
$$g_2=x-\frac{9}{8}y-\frac{3}{8}z-\frac{1}{2}$$

\section{Conclusion}
In \cite{zhang_2}, a bridge overriding the gap between approximate
computation and exact results was built. However, the algorithm in
\cite{zhang_2} requires the error between approximation and exact
result is less than $1/((2N+2)N(N-1))$. In this paper, we propose
an algorithm that only requires the error is less than
$1/(4N(N-1))$, which decreases the cost of computing
approximation.  Just like the algorithm in \cite{zhang_2}, our
method can be applied in many aspect, such as proving inequality
statements and equality statements, and computing resultants, etc.
Thus we can take fully advantage of approximate methods to solve
larger scale symbolic computation problems.

\end{document}